\documentclass[12pt]{article} 

\usepackage{amsmath, amsfonts, amssymb}
\usepackage{mathrsfs}
\usepackage{theorem}
\usepackage{bm}
\usepackage[dvipdf]{graphicx}

\RequirePackage[colorlinks,citecolor=blue,urlcolor=blue,bookmarksopen]{hyperref}

\newtheorem{mythm}{Theorem}[section]
\newtheorem{myprop}[mythm]{Proposition}
\newtheorem{mylem}[mythm]{Lemma}
\newtheorem{mycor}[mythm]{Corollary}

{\newtheorem{myrem}[mythm]{Remark}}
{\newtheorem{myexam}{Example}[section]}

\newcommand{\dis}{\displaystyle}

\def\R{\mathbb R}
\def\Z{\mathbb Z}
\def\N{\mathbb N}

\def\B{\mathscr B}

\def\F{\mathscr F}
\def\d{\text{\rm{d}}}
\def\E{\mathbb E}
\def\p{\mathbb P}
\def\e{\text{\rm{e}}}

\def\la{\langle}
\def\raa{\rangle}
\def\La{\Lambda}
\def\veps{\varepsilon}

\def\S{\mathcal S}

\def\pb{\mathscr{P}}

\def\wt{\widetilde}

\def\var{\mathrm{var}}

\newcommand\ea{{\veps,\alpha}}

\newcommand{\fin}{\hfill $\square$\par}

\numberwithin{equation}{section}

\allowdisplaybreaks

\begin{document}

\title{Averaging principle for two time-scale regime-switching
processes\footnote{Supported in
 part by NNSFs of China (No. 12271397,  11831014)  and National Key R\&D Program of China (No. 2022YFA1000033)}}

\author{Yong-Hua Mao\thanks{School of Mathematical Sciences, Beijing Normal University, Beijing, China} and Jinghai Shao\thanks{Center for Applied Mathematics, Tianjin University, Tianjin, China.   Email: shaojh@tju.edu.cn.}
}
\maketitle

\begin{abstract}
This work studies the averaging principle for  a fully coupled two time-scale system, whose slow process is a diffusion process and fast process is a purely jumping process on an infinitely countable state space.
The ergodicity of the fast process has important impact on the limit system and the averaging principle. We showed that under strongly ergodic condition, the limit system admits a unique solution, and the slow process converges in the $L^1$-norm to the limit system. However, under certain weaker ergodicity condition, the limit system admits a solution, but not necessarily unique, and the slow process can be proved to converge  weakly to a solution of the limit system.
\end{abstract}



\section{Introduction}
We study  in this work a fully coupled two time-scale stochastic system $(X_t^\ea, Y_t^\ea) $ in  $\R^d\times\S$, where $\S=\{1,2,\ldots,N\}$ with $N\leq \infty$. The slow process $(X_t^\ea) $ is  described as a solution to the following stochastic differential equation (SDE):
\begin{equation}\label{o-1}
\begin{split}
  \d X_t^{\veps,\alpha}&=b(X_t^{\veps,\alpha}, Y_t^{\veps, \alpha})\d t+\sqrt{\veps}\sigma(X_t^{\veps,\alpha}, Y_t^{\veps, \alpha})\d W_t,\\
  X_0^{\veps,\alpha}&=x_0\in\R^d,\quad Y_0^{\veps,\alpha}=i_0\in\S,
\end{split}
\end{equation}
and the fast process $(Y_t^{\veps,\alpha})$ is a jumping-process on $\S$ satisfying
\begin{equation}\label{o-2}
\p(Y_{t+\delta}^{\veps,\alpha}=j|Y_t^{\veps,\alpha} =i,X_t^{\veps,\alpha}=x)=\begin{cases} \frac 1\alpha q_{ij}(x)\delta+o(\delta),\ &\text{if $i\neq j$},\\
1+\frac 1\alpha q_{ii}(x)\delta+o(\delta),&\text{if $i=j$}
\end{cases}
\end{equation}
for $\delta>0$, $i,\,j\in\S$, $x\in \R^d$, and $\veps,\,\alpha$ are small positive parameters.  In the existing literatures, the system $(X_t^\ea,Y_t^\ea)$ is called \emph{fully coupled} if the diffusion coefficient $\sigma$ of slow process $(X_t^\ea)$ depends on the fast process $(Y_t^\ea)$ and the transition rates $(q_{ij}(x))_{i,j\in\S}$ of the fast process $(Y_t^\ea)$   depends on $(X_t^\ea)$ as well.

Multi-scale systems arise in many research fields such as in  biology systems \cite{FGC,KK,KKP,KhY,Pop}, in mathematical  finance  \cite{FFF,FFK}, etc. Correspondingly, there are many works devoted to the study of averaging principle, central limit theorems, and large deviations of these stochastic models. For a two time-scale system where both slow and fast components are continuous processes given as solutions of SDEs, these problems have been extensively studied, such as,  in \cite{Bao,Kh68,KhY,Kh05,Lip,Liu,Puh,Ver1,Ver2}, in \cite{HL20} for SDEs driven by fractional Brownian motions. The interaction between the fast  component and the slow one makes a fully coupled two time-scale system much complicated, which has been  revealed in the works \cite{Kifer,Puh,Ver1,Ver2}.

The averaging principle says that the slow process $(X_t^\ea)$ will converge to some limit process $(\bar{X}_t)$ as $\veps,\alpha\to 0$. When the fast process $(Y_t^\ea)$ does not depend on $(X_t^\ea)$, usually called an uncoupled system, the averaging principle often holds in quite general conditions. However, when $(Y_t^\ea)$ depends on $(X_t^\ea)$ and particularly $(Y_t^\ea)$ does not locate in a compact space, it becomes more difficult to establish the averaging principle. In this work we focus on addressing the impact on the limit behavior of $(X_t^\ea,Y_t^\ea)$ caused by: 1)  various ergodicities of the fast process on the wellposedness of the limit process $(\bar{X}_t)$; 2) when the state space $\S$ is infinitely countable, the dependence on the fixed state of the slow process $(X_t^\ea)$ of the invariant measure of $(Y_t^\ea)$.

Let us review some known works in the setting similar to ours.
In the situation that $(Y_t^\ea)$ is a continuous time Markov chain independent of the slow process $(X_t^\ea)$, Eizenberg and Freidlin \cite{EF93}, Freidlin and Lee \cite{FL96}  investigated separately the limit behavior of solutions of PDE systems with Dirichlet boundary associated with $(X_t^\ea,Y_t^\ea)_{t\geq 0}$ when the  diffusion coefficient of $X_t^\ea$ does not depend or depends  on $Y_t^\ea$. These two works reveal that whether the diffusion coefficient of $X_t^\ea$ depends on $Y_t^\ea$ or not has important impact on the method to study the limit behavior of $(X_t^\ea,Y_t^\ea)$.  To provide a decisive estimate on the difference between $(X_t^\ea)$ and its limit process, a large deviation principle (LDP) was established in
\cite{HYZ,HMS}.

For a setting where the fast process $(Y_t^\ea)$ is a jumping process depending on the slow process $(X_t^\ea)$ as well,
the averaging principle and LDP have been studied by Faggionato, Gabrielli, and Crivellari \cite{FGC} and Budhiraja, Dupuis and Ganguly \cite{BDG18}. \cite{FGC} considered  a simple case without diffusion term for the slow component by the nonlinear semigroup method developed by Feng and Kurtz \cite{FK06}. Whereas, \cite{BDG18} considered a fully coupled case by using the weak convergence method, and established a process level large deviation  principle.

All the aforementioned works, no matter whether the fast jumping process $(Y_t^\ea)$ depends on $(X_t^\ea)$ or not, considered only the situation that the state space $\S$ of  $(Y_t^\ea)$ is a finite state space, which is hence compact.  However, the infinite countability of the state space $\S$ of $(Y_t^\ea)$ has important impact on the averaging principle and LDP of $(X_t^\ea,Y_t^\ea)$.
For example, in a simple setting $\alpha\equiv 1$,  Bezuidenhout \cite{Bez} studied the LDP of certain  functionals of $(X_t^\ea, Y_t^\ea)$ with the diffusion coefficient of $(X_t^\ea)$ independent of $(Y_t^\ea)$. It showed that  the LDP holds when   $(Y_t^\ea)$ is in a finite state space. Furthermore, it was shown by a counterexample that when $(Y_t^\ea)$ is a Markov chain in an infinite state space, the LDP may fail. Meanwhile,
as our studied system $(X_t^\ea,Y_t^\ea)$ is fully coupled, the invariant probability measure $\pi^x=(\pi^x_i)_{i\in\S}$ of $(Y_t^\ea)$ will depend on the position $x$ of the slow process $(X_t^\ea)$. The infinite countability of $\S$ makes  the regularity of $x\mapsto \pi^x$ become much more complicated than the case that $\S$ is finite. The regularity of $x\mapsto \pi^x$ has important impact on the characterization of the limit system.

Precisely, suppose $(q_{ij}(x))_{i,j\in\S}$ is a conservative, irreducible transition rate matrix for every $x\in\R^d$, which is Lipschitz continuous in $x$ in certain matrix norm. Let $\pb(\S)$ be the space of all probability measures over $\S$ endowed with the total variation norm. Let $\pi^x\in \pb(\S)$ be the invariant probability measure associated with $(q_{ij}(x))_{i,j\in\S}$ provided it  exists. Then, when $\S$ is a finite state space, $x\mapsto \pi^x$ as a function from $\R^d$ to $\pb(\S)$ is Lipschitz continuous. This result has been  proved in \cite{FGC} and \cite{BDG18} in different ways.  \cite{FGC} proved it by the Perron-Frobenius theorem to express  $\pi^x$ in terms of a nonzero right eigenvector of $(q_{ij}(x))_{i,j\in\S}$ corresponding to the eigenvalue 0. In \cite{BDG18}, it is proved through expressing $\pi^x$ as a polynomial of transition probabilities according to Freidlin and Wentzell \cite{FW}.  Nevertheless, these two methods are infeasible when $\S$ is infinite. Moreover, when $\S$ is infinitely countable, $x\mapsto \pi^x$ could be not Lipschitz continuous and even not H\"older continuous of any exponent in $(0,1)$; see our Example \ref{exam-1} below.

To establish the averaging principle when $\S$ is infinitely countable, our main challenge is to  study the regularity of $x\mapsto \pi^x$ from $\R^d$ to $\pb(\S)$. To overcome this difficulty, the ergodic property of $P_t^x$  plays a crucial role, where $P_t^x$ denotes the semigroup associated with the Markov chain with transition rate matrix $(q_{ij}(x))_{i,j\in\S}$.
We shall show that $x\mapsto \pi^x$ is Lipschitz continuous if $P_t^x$ is \emph{strongly ergodic}   uniformly w.r.t. $x$ based on an integration by parts formula for continuous time Markov chains. If supposing only that  $P_t^x$ is \emph{ergodic} and  $\|P_t^x(i,\cdot)-\pi^x\|_\var\leq C_i\eta_t$ for $i\in \S$ with $C_i>0$, $\eta_t\in [0,2]$ satisfying $\int_0^\infty \eta_s\d s<\infty$,   $x\mapsto \pi^x$ is shown to be $1/2$-H\"older continuous. To prove this assertion,  we develop a coupling method for parameter-dependent Markov chains based on Skorokhod's representation theorem for jumping processes.
Consequently,
under the strongly ergodic condition, the equation to characterize the limit process $(\bar{X}_t)$ admits a unique solution, and we can show that $(X_t^\ea)$ converges in $L^1$-norm   to $(\bar{X}_t)$ as $\veps,\alpha \to 0$. However, under ergodic  condition,   $(X_t^\ea)$   converges weakly to its limit process provided that the limit system is unique. The ratio $\veps/\alpha$ as $\veps,\alpha\to 0$ has no impact on the averaging principle. Nevertheless, the large deviation principle of $(X_t^\ea,Y_t^\ea)$ will be shown to depend heavily on the ratio $\veps/\alpha$ in our another work.

The remainder of this work is organized as follows. In Section 2, we state the main results of this work including: the regularity of $x\mapsto \pi^x$ under two different ergodicity conditions, and the averaging principle for $(X_t^\ea, Y_t^\ea)_{t\geq 0}$ in respectively strong and weak convergence sense.
Section 3 is devoted to developing the coupling method for parameter-dependent Markov chains, which is not only the basis to study the regularity of $x\mapsto \pi^x$ under the weak ergodicity condition of $(Y_t^\ea)$, but also plays an important role to decouple the close interaction between $(X_t^\ea)$ and $(Y_t^\ea)$ to establish the averaging principle.
The arguments of main results are all presented in Section 4.




\section{Statement of main results}\label{sec-2}

This section is devoted to establishing the averaging principle for $(X_t^\ea,Y_t^\ea)_{t\geq 0}$ as $\veps,\,\alpha$ go to zero.
%
Let us begin with introducing three fundamental conditions  on the stochastic system $(X_t^\ea,Y_t^\ea)$, which  will be used throughout this work.
\begin{itemize}
  \item[$\mathrm{(A1)}$] There exist constants $K_1,\,K_2>0$ such that
      \begin{gather*}
        |b(x,i)-b(y,i)| +\|\sigma(x,i)-\sigma(y,i)\|\leq K_1|x-y|,\\
        |b(x,i)|+\|\sigma(x,i)\| \leq K_2,\quad \ \ \ x,y\in \R^d,\ i\in \S.
      \end{gather*}
  \item[$\mathrm{(A2)}$] For each $x\in \R^d$, $(q_{ij}(x))_{i,j\in\S}$ is a conservative, irreducible transition rate matrix. Assume $\kappa:=\sup_{i\in\S} \sum_{j\in\S,j\neq i}\sup_{x\in \R^d} q_{ij}(x)<\infty$.
  \item[$\mathrm{(A3)}$] There exists a constant $K_3>0$  such that
      \begin{gather*}
      \|Q(x)-Q(y)\|_{\ell_1}:=\sup_{i\in\S}\sum_{ j\neq i}|q_{ij}(x)-q_{ij}(y)|\leq K_3|x-y|,\quad x,y\in \R^d.
      \end{gather*}
\end{itemize}
Under these conditions (A1)-(A3), the two time-scale system \eqref{o-1}, \eqref{o-2} admit a unique strong solution to any initial value $X_0^\ea=x_0\in \R^d$ and $Y_0^\ea=i_0\in \S$; see, e.g. \cite{YZ} or \cite{Sh15} under certain more general non-Lipschitz conditions. To focus our idea on the impact of the ergodicity of $(Y_t^\ea)$ on the averaging principle,  we impose a simple condition (A1) on the slow process $(X_t^\ea)$. We refer the readers to \cite{Liu} for the technique  to generalize (A1) to the local Lipschitz condition.

%

For the fully coupled two time-scale system  $(X_t^\ea,Y_t^\ea)$, in contrast to uncoupled two time-scale systems, the regularity of invariant probability measure $\pi^x$ associated with the $Q$-matrix $(q_{ij}(x))_{i,j\in\S}$ increases the complexity and difficulty of characterizing the limit system $(\bar{X}_t)$ of $(X_t^\ea,Y_t^\ea)$ as $\veps,\,\alpha\to 0$.
As mentioned in the introduction, when $\S$  is a finite state space, and $(q_{ij}(x))_{i,j\in\S}$ is Lipschitz continuous in $x$, then its associated invariant probability measure $\pi^x=(\pi_i^x)_{i\in\S}$ is also Lipschitz continuous in $x$, which has been proved in \cite{BDG18,FGC}. However, when $\S$ is infinitely countable, this becomes uncertain. Note that the invariant probability measure is also a left eigenvector to the $Q$-matrix. The perturbation on  linear generators can cause significant changes on its corresponding eigenvalues and eigenvectors. To see the complexity of this problem, one can refer to the fruitful researches on the perturbation theory of linear operators; see, for instance, the monograph \cite{Kato} and references therein.

Let us recall some notations on the ergodicity of Markov chains (cf. \cite{Chen,MT}). Let $P_t $ denote a semigroup associated with a continuous time Markov chain on the state space $\S$. Suppose that there exists an invariant probability measure $\pi=(\pi_i)_{i\in\S}$. The total variation distance between $P_t(i,\cdot)$ and $\pi$ is defined by
\[\|P_t(i,\cdot)\!-\!\pi\|_\var\!=2\sup\big\{P_t(i, A)\!-\!\pi(A);\, A\!\in \! \mathscr{B}(\S)\big\}=\sup\big\{|P_t(i,f)\!-\!\pi(f)|; |f|\leq 1\big\},
\] where $\mu(f):=\sum_{i\in\S}\mu_if(i)$ for any probability measure $\mu$ on $\S$. The Markov chain is called \emph{ergodic} if
\[\lim_{t\to \infty} \|P_t(i,\cdot)-\pi\|_{\var}=0, \ \quad  i\in\S;\]
the process is called \emph{exponentially ergodic}, if
\[ \|P_t(i, \cdot)-\pi\|_{\var}\leq C_i \e^{-\veps t}, \ t>0, \ \text{for some $\veps>0$, constants $C_i>0$, $i\in\S$;}
\]
the process  is called \emph{strongly ergodic} or \emph{uniformly ergodic}, if
\[\lim_{t\to \infty} \sup_{i\in\S} \|P_t(i,\cdot)-\pi\|_\var=0.\]
It is known that if the chain is strongly ergodic, its convergence rate must be of exponential type, i.e.
\[ \sup_{i\in\S}\|P_t(i,\cdot)-\pi\|_{\var} \leq C\e^{-\lambda t}, \quad t>0,\]
for some constants $C,\lambda>0$; see, for example,  \cite[Lemma 4.1]{Mao02}.
Consequently, it is easy to see that ergodic Markov chain on a finite state space must be strongly ergodic. Accordingly, we first generalize the results in \cite{BDG18,FGC} for Markov chains on a finite state space to the setting on an infinite state space under the strongly  ergodic  condition.

Let $P_t^x$ be the semigroup associated with the $Q$-matrix $(q_{ij}(x))_{i,j\in \S}$, and $\pi^x$ its associated invariant probability measure provided it exists throughout this work.

\begin{itemize}
  \item[$\mathrm{(A4)}$]  Suppose that $P_t^x$ is strongly ergodic uniformly in $x$, that is,  there exist constants $c_1,\lambda_1>0$ such that
      \[ \sup_{i\in\S} \|P_t^x(i,\cdot) -\pi^x\|_{\var}\leq c_1\e^{-\lambda_1 t},\quad \forall\,t>0,\ x\in \R^d.\]
\end{itemize}

\begin{myprop}[Strongly ergodic case]\label{prop-1}
Assume (A2), (A3) and (A4) hold.  Then, the functional $\R^d\ni x\mapsto \pi^x\in \pb(\S)$   is Lipschitz continuous, i.e.
\[\|\pi^x-\pi^y\|_{\var}\leq C_\pi|x-y|,\quad\quad  x,y\in \R^d,\]
where $C_\pi=\frac{4c_1K_3}{\lambda_1}$ and constants $c_1,\lambda_1$ given in (A4), $K_3$ given in (A3).
\end{myprop}
To make the presentation transparent, we defer the argument to Section \ref{argument}. It is useful to mention the works \cite{Chen,Mao02} and references therein, which provide various sufficient conditions for  strong ergodicity of continuous-time Markov chains and diffusion processes.

We proceed to investigate the regularity of $x\mapsto \pi^x$ under certain ergodic condition weaker than strong ergodicity. Unfortunately, under weaker ergodic condition and without the uniformity w.r.t.\! $x$, the Lipschitz continuity of $x\mapsto \pi^x$ in the total variation norm may fail. To illustrate it, we construct  an example as follows.

\begin{myexam}\label{exam-1}
For each $x\in (0,1)$, let $(Y_t^x)_{t\geq 0}$ be a birth-death process on $\S=\{1,2,\ldots\}$ with birth rate $q_{ii+1}(x)=b_i(x) =x$ for $i\geq 1$ and death rate $q_{ii-1}(x)=a_i(x)=1$ for $i\geq 2$. It is clear that $q_{ij}(x)$ is Lipschitz continuous in $x$ for all $i,j\in\S$. Then,
\begin{itemize}
  \item[$\mathrm{(i)}$] for each $x\in (0,1)$, the birth-death Markov chain $(Y_t^x)_{t\ge 0}$ is exponentially ergodic, but not strongly ergodic, satisfying
      \begin{equation}\label{e-exam-1}
      \|P_t^x(i,\cdot)-\pi^x\|_\var\leq C_i(x)\e^{-( 1-\sqrt{x})^2t},\quad t>0, \, i\in \S,
      \end{equation}
      for  some positive constants $C_i(x)$ depending on $i\in\S$ and $x\in (0,1)$.
  \item[$\mathrm{(ii)}$] Its invariant probability measure $\pi^x=(\pi^x_i)_{i\geq 1}$ is given by
      \begin{equation}\label{e-exam-2}
      \pi_i^x= (1-x)x^{i-1},\quad i\geq 1,
      \end{equation}
      and for any $\beta\in (0,1]$
      \begin{equation}\label{e-exam-3}
      \sup_{x\neq y}\frac{\|\pi^x-\pi^y\|_{\var}}{|x-y|^\beta}=\infty.
      \end{equation}
      This means that  $x\mapsto \pi^x$ is \textbf{not} H\"older continuous of any exponent $\beta\in (0,1]$.
\end{itemize}
\end{myexam}
The argument of assertions stated in Example \ref{exam-1} is also deferred to Section \ref{argument}.

Now, let us consider the following ergodic condition weaker than strong ergodicity condition (A4).
\begin{itemize}
  \item[$\mathrm{(A5)}$]  Assume that there exist   a  positive  function $\theta:\S\to (0,\infty)$, a decreasing function $\eta:[0,\infty)\to [0,2]$ satisfying $ \int_0^\infty \eta_s\d s<\infty$ such that    \[ \|P_t^x(i ,\cdot)-\pi^x\|_{\var}\leq \theta(i)\eta_t,\quad \ t\geq 0, \ x\in\R^d,\  i\in\S.\]
\end{itemize}

%

\begin{myprop}\label{prop-2}
Assume the conditions (A2), (A3) and (A5) hold, then $x\mapsto \pi^x$ is $1/2$-H\"older  continuous, i.e.
\begin{equation}\label{aa-5}
\|\pi^x-\pi^y\|_{\var} \leq K_4\sqrt{|x-y|},\qquad x,y\in\R^d,
\end{equation} where $K_4=\sqrt{K_3(\inf_{i\in\S}\theta(i)) \int_0^\infty \!\!\eta_s\d s } $.
\end{myprop}

This proposition is proved based on an intricate construction of  coupling process of $(\wt{Y}_t^x)  $ and $(\wt{Y}_t^y) $ with $Q$-matrix $(q_{ij}(x))_{i,j\in\S}$ and $(q_{ij}(y))_{i,j\in\S}$ respectively in terms of Skorokhod's representation for jumping processes, which is presented in Section \ref{coupling}. Our construction method in current work improves the one used in \cite{SY} to study the  stability of regime-switching processes under the perturbation of $Q$-matrix and in \cite{SZ21} to study the continuous dependence of intial values for stochastic functional differential equations with state-dependent regime-swtiching. The key point is the estimate of
$\frac 1t\int_0^t\p(\wt{Y}_s^x\neq \wt{Y}_s^y)\d s$ in terms of the difference between $(q_{ij}(x))_{i,j\in\S}$ and $(q_{ij}(y))_{i,j\in\S}$.

Next, we go to establish the averaging principle for  $(X_t^\ea, Y_t^\ea)$ as $\veps,\alpha\to 0$. Let
\begin{equation}\label{a-3}
\bar{b}(x)=\sum\nolimits_{i\in\S} b(x,i)\pi_i^x,
\end{equation}
and the limit system of $(X_t^\ea,Y_t^\ea)$ will be given as the solution to the ordinary differential equation (ODE)
\begin{equation}\label{a-4}
\d \bar{X}_t=\bar{b}(\bar{X}_t)\d t,\quad\quad \bar{X}_0=x_0.
\end{equation}
Under conditions (A1) and (A4), by Proposition \ref{prop-1}, it is easy to see $\bar{b}$ is Lipschitz continuous, and hence ODE \eqref{a-4} admits a unique solution. Under the strongly ergodic condition (A4), we can get  $L^1$-convergence of $X_t^\ea$ to $\bar{X}_t$ as $\veps, \alpha\to 0$.

\begin{mythm}\label{thm-1}
Assume (A1)-(A4) hold. Let $(X_t^\ea, Y_t^\ea)$ be the solution to \eqref{o-1}, \eqref{o-2}, and $(\bar{X}_t)$ the solution to \eqref{a-4}.  Then
\begin{equation*}
  \lim_{\veps,\alpha\to 0} \E|X_t^\ea-\bar{X}_t|=0,\quad t>0.
\end{equation*}
\end{mythm}
However, under (A1) and (A5), by Proposition \ref{prop-2}, $\bar{b}$ can be shown only to be H\"older continuous just as   $\pi^x$.  In this situation, thanks to Peano's theorem, ODE \eqref{a-4} admits a solution, but may loss the uniqueness. Consequently, under the weaker ergodic condition (A5) the limit system $(\bar{X}_t)$ becomes more complicated, and $(X_t^\ea)$ can be shown to converge weakly to its limit whenever ODE admits a unique solution. The precise result is given in the following theorem.

\begin{mythm}\label{thm-1.5}
Assume that (A1)-(A3) and (A5) hold. In addition, when $\S$ is infinitely countable, suppose that there exist constants $c_2>0,\,c_3<\infty$ such that the function $\theta(\cdot)$ given in (A5) also satisfies
\begin{equation}\label{cond-1}
Q(x)\theta(i)=\sum_{j\in\S}q_{ij}(x)\theta(j)
      \leq -c_2\theta(i)+c_3,\quad x\in \R^d, \ i\in\S.
\end{equation}
Let $(X_t^\ea, Y_t^\ea)$ be the solution to \eqref{o-1}, \eqref{o-2}.
Then, for each $T>0$, the  set of distributions of $\{(X_t^\ea)_{t\in [0,T]};\veps,\alpha\in (0,1)\}$ in $\mathcal{C}([0,T];\R^d)$ is tight, and any convergent subsequence of $\{(X_t^\ea)_{t\in [0,T]};\veps,\alpha>0\}$ shall converge weakly to a solution $(\bar{X}_t)_{t\in [0,T]}$ of ODE \eqref{a-4}. Moreover, if ODE \eqref{a-4} admits a unique solution, then $(X_t^\ea)_{t\in [0,T]}$ converges weakly to the unique solution $(\bar{X}_t)_{t\in [0,T]}$ of ODE \eqref{a-4} as $\veps,\alpha\to 0$.
\end{mythm}


\begin{myrem}
  Theorems \ref{thm-1} and \ref{thm-1.5} tell  us a fundamental fact: the limit system $(\bar X_t)$ and   the convergence of $(X_t^\ea)$ to this limit system  do  not depend on the ratio  $\veps /\alpha$ as $\veps,\alpha\to 0$.
\end{myrem}

\begin{myexam}\label{exam-3}
\begin{enumerate}
  \item Let  $(Y_t^x)$ is a continuous-time Markov chain on $\S=\{1,2,\ldots\}$ with the transition rate matrix
      \begin{equation*}
      \begin{split}
        q_{ij}(x)&=\big(1-\e^{-|x|-\alpha}\big)\e^{-(j-1)(|x|+\alpha)},\quad j\neq i,\\
        q_{ii}(x)&=-\big(1-\big(1-\e^{-|x|-\alpha}\big)\e^{ -(i-1)(|x| +\alpha)}\big),
      \end{split}
      \end{equation*} where $\alpha>0$ for $x\in \R$. Then, according to \cite[Theorem 4.45]{Chen-1}, the process $(Y_t^x)$ satisfies the condition (A4). Moreover,
      \begin{align*}
        \|Q(x)-Q(y)\|_{\ell_1}&:=\sup_{i\geq 1}\sum_{j\neq i} |q_{ij}(x)-q_{ij}(y)|
         \leq 2\big(\sum_{j\geq 1}j\e^{-j\alpha}\big) |x-y|, \ x,y\in\R,
      \end{align*} which means that (A3)   holds as well.
   \item Let $(Y_t^x)$ be associated with the transition rate matrix $(q_{ij}(x))$ given by
       \[q_{i(i+1)}(x)=2+\sin x,\  q_{i1}(x)=2-\sin x, \  q_{ii}(x)=-4; \ q_{ij}(x)=0, \text{otherwise}\] for $ i\geq 1, \ x\in \R$. Again, by \cite[Theorem 4.45]{Chen-1}, $(Y_t^x)$ satisfies (A4).
   \item Let $(Y_t^x)$ be a birth-death process on $\S=\{1,2,\ldots\}$ with $b_i(x)=q_{ii+1}(x)=1$ for $i\geq 1$, $a_i(x)=q_{ii-1}(x)=2-\frac 12\sin x$ for $i\geq 2$, $x\in \R$. Then $(Y_t^x)$ is exponentially ergodic and satisfies (A5).
\end{enumerate}
\end{myexam}

\section{Construction of the coupling processes}\label{coupling} In this part we introduce the coupling processes used in the study of regularity of $x\mapsto \pi^x$ and in decoupling the interaction between the slow process $(X_t^\ea)$ and the fast process $(Y_t^\ea)$ in order to establish the averaging principle. This part deals with the technical difficulties caused by the full dependence between $(X_t^\ea)$ and $(Y_t^\ea)$. In the spirit of Skorokhod, we express a state-dependent jumping process over $\S$ in terms of an integral w.r.t.\! a Poisson random measure. In order to deal with the case $\S$ being infinitely countable, we modify the construction method of intervals used in Skorokhod's representation theorem, which is quite different to the extensively used one (cf. e.g. \cite{Gho,Sh15,YZ}).

Consider the  solutions $(  X_t^x,Y_t^x)$ and $(\wt X_t^y, \wt Y_t^y)$ respectively to the following SDEs:
\begin{equation}\label{a-4.8}
\begin{aligned}
  &\d X_t^x=b(  X_t^x, Y_t^x)\d t+\sigma( X_t^x, Y_t^x)\d W_t,\quad X_0^x=x\in \R^d, Y_0^x=i_0\in \S,\\
  &\p(Y_{t+\delta}^x=j|Y_t^x=i, X_t^x=z)=\begin{cases}
     {q}_{ij}(z)\delta+o(\delta),\ &i\neq j,\\
    1+  q_{ii}(z)\delta+o(\delta), &i=j,
  \end{cases}
\end{aligned}
\end{equation}
and
\begin{equation}\label{a-4.9}
\begin{aligned}
  &\d \wt{X}_t^y =f(\wt{X}_t^y, \wt{Y}_t^y)\d t+g(\wt{X}^y_t,\wt{Y}_t^y)\d W_t,\ \ \wt{X}_0^y=y\in \R^d, \wt{Y}_0^y=i_0\in\S,\\
  &\p(\wt{Y}_{t+\delta}^y=j|\wt{Y}_t^y=i, \wt{X}_t^y=z)=\begin{cases}
     {q}_{ij}(z)\delta+o(\delta),\ &i\neq j,\\
    1+  q_{ii}(z)\delta+o(\delta), &i=j,
  \end{cases}
\end{aligned}
\end{equation} for $\delta>0$.

\begin{mylem}[Key lemma]\label{lem-1}
Suppose that (A1), (A2) hold  and $f,g$  satisfy (A1) replacing $b$ and $\sigma$ respectively. For every $x,y\in \R^d$, $x\neq y$ and every $i_0\in \S$, there is a coupling process $(X_t^x, Y_t^x)_{t\geq 0}$ and $(\wt{X}^y_t, \wt Y_t^y)_{t\geq 0}$  satisfying SDEs \eqref{a-4.8} and \eqref{a-4.9} respectively such that
\begin{equation}\label{a-5}
\frac 1t\int_0^t\p(  Y_s^x\neq \wt Y_s^y)\d s\leq \!
  \int_0^t\!\E\big[\|Q(X_s^x)\!-\!Q(\wt{X}_s^y) \|_{\ell_1}  \big] \d s,\quad t>0,
\end{equation}
where $\|Q(x)\!-\!Q(y)\|_{\ell_1}\!=\!\sup\limits_{i\in\S}\! \sum\limits_{j\in\S,j\neq i} |q_{ij}(x)\!-\!q_{ij}(y)|$.
\end{mylem}

As an application of Lemma \ref{lem-1}, consider a special case:
$b=f=0$, $\sigma=g=0$, then $X_t^x\equiv x$, $\wt{X}_t^y\equiv y$, and we obtain that:

\begin{mycor}\label{cor-1}
Under (A2), for every $x,y\in \R^d$,  there is a coupling process $(Y_t^x, \wt{Y}_t^y)$ associated respectively with the $Q$-matrix $(q_{ij}(x))_{i,j\in\S}$ and $(q_{ij}(y))_{i,j\in\S}$ such that
 \begin{equation}\label{a-5.1}
\frac1t \int_0^t\!\p(Y_s^x\neq \wt{Y}_s^y)\d s\leq   t  \|Q(x)-Q(y)\|_{\ell_1} , \quad t>0.
\end{equation}
\end{mycor}

Corollary \ref{cor-1} tells us that when $|x-y|$ tends to 0, we can construct a coupling process $(Y_t^x,Y_t^y)$ such that $\frac1t \int_0^t\!\p(Y_s^x\neq \wt{Y}_s^y)\d s$ goes to $0$ when $x\mapsto Q(x)$ is continuous.

\noindent \textbf{Argument of Lemma \ref{lem-1}}  We need first construct suitable intervals related to the transition rate matrix $(q_{ij}(x))_{i,j\in\S}$ so as to express the jumping processes $(Y_t^x)$ and $(\wt{Y}_t^y)$ in terms of a common Poisson random measure. The proof is divided into two steps.

\noindent\textbf{Step 1}.
The first step is to construct a sequence of intervals associated with the transition rate matrix $(q_{ij}(z))$, $z\in \R^d$. Our construction method is applicable when $\S$ is finite or infinite, and  is more suitable to cope with the case $\S$ is infinite than the construction method used in \cite{Gho,Sh15,YZ}.

  Precisely, let $\gamma_n=\sup_{k\neq n}\sup_{z\in \R^d} q_{nk}(z)$ for $n\in \S$, and by (A2) we get $\gamma_n\leq \kappa<\infty$ for all $n\in\S$.
  Let $\displaystyle \Gamma_{1k}(z) \!=\![(k-2)\gamma_1,(k-2)\gamma_1 +q_{1k}(z))$ for $k\geq 2$, and for $n\geq 2$,
  \begin{align*}
    \Gamma_{nk}(z)&\!=\![ (k-n-1)\gamma_n, (k -n-1)\gamma_n +q_{nk}(z)),\ \text{if}\,  k>n,\\
    \Gamma_{nk}(z)&\!=\![(k+1-n)\gamma_n-q_{nk}(z), (k+1-n)\gamma_n),\ \text{if}\ 1\leq k<n,
  \end{align*}
  and
  \begin{equation}\label{U}
    U_n=\bigcup_{z\in \R^d}\bigcup_{k\geq 1,k\neq n}\Gamma_{nk}(z),\quad n\geq 1,
  \end{equation}
  where   $\kappa$ is given in   (A2). For notation convenience, we put $\Gamma_{ii}(z)=\emptyset$ and $\Gamma_{ij}(z)=\emptyset$ if $q_{ij}(z)=0$, $i,j\in\S$, $z\in\R^d$. Due to (A2),
  \begin{equation}\label{U-1}
  \mathbf{m}(U_n)\leq \sum_{k\geq 1, k\neq n}\sup_{z\in \R^d}q_{nk}(z)\leq \kappa,
  \end{equation} where   $\mathbf{m}(\d x)$ denotes the Lebesgue measure over $\R$.

  Secondly, we provide an explicit construction of the Poisson random measure  as in \cite{IW}, which helps us to illustrate the calculation below.
   Let
  $\xi_i^{(k)}$, $k,i=1,2,\ldots$, be $U_k$-valued random variables with
  \[\p(\xi_i^{(k)}\in \d x)=\frac{\mathbf{m}(d x)}{\mathbf{m}(U_k)},\]
  and $\tau_i^{(k)}$, $k,i\geq 1$, be non-negative random variables satisfying $\p(\tau_i^{(k)}>t)=\exp[-t\mathbf{m}(U_k)]$, $t\geq 0$. Suppose that $\{\xi_i^{(k)},\tau_i^{(k)}\}_{i,k\geq 1}$ are all mutually independent. Put
  \[\zeta_n^{(k)}=\tau_1^{(k)}+\cdots+\tau_n^{(k)} \ \text{for}\ n,\,k\geq 1, \ \text{and}\ \zeta_0^{(k)}=0,\ k\geq 1.\]
  Let
  \[D_{\mathbf{p}}= \bigcup_{k\geq 1}\bigcup_{n\geq 0}\big\{\zeta_n^{(k)}\big\},\]
  and
  \[ \mathbf{p}(t)=\sum_{0\leq s< t} \Delta \mathbf{p}(s), \quad \Delta\mathbf{p}(s)=0\ \text{for $s\not\in D_{\mathbf{p}}$},\  \Delta\mathbf{p}(\zeta_n^{(k)})  =\xi_n^{(k)},\quad k,n\geq 1,\]
  where $\Delta\mathbf{p}(s)=\mathbf{p}(s)- \mathbf{p}(s-)$.
  Correspondingly, put
  \[\mathcal{N}_{\mathbf{p}}([0,t]\times A)=\#\{s\in D_{\mathbf{p}};\ 0<s\leq t,\Delta\mathbf{p}(s)\in A\},\ t>0,\, A\in\mathscr{B}([0,\infty)).\]
  As a consequence, we get a Poisson point process $(\mathbf{p}(t))$ and a Poisson random measure $\mathcal{N}_{\mathbf{p}}(\d t,\d x)$ with intensity $\d t\mathbf{m}(\d x)$.

  Thirdly, let  \[\vartheta(x,i,z)=\sum_{j\in\S}(j-i) \mathbf{1}_{\Gamma_{ij}(x)}(z).\]
  The desired coupling process is defined as the solutions to the following SDEs.
  \begin{equation}\label{ac-1}
  \begin{cases}
  \d X_t^x&=b(X_t^x, Y_t^x)\d t+\sigma(X_t^x, Y_t^x)\d W_t,\\
    \d  Y_t^x&=\int_{[0,\infty)}\! \! \vartheta(X_t^x,Y_{t-}^x,z)
    \mathcal{N}_{\mathbf{p}}(\d t,\d z),\quad X_0^x=x,\, Y_0^x=i_0.
  \end{cases}
  \end{equation}
  \begin{equation}\label{ac-2}
  \begin{cases}
    \d \wt{X}_t^y&=f(\wt{X}_t^y, \wt{Y}_t^y)\d t+g(\wt{X}_t^y, \wt{Y}_t^y)\d W_t,\\
    \d \wt{Y}_t^y&=\int_{[0,\infty)}\!\!\vartheta( \wt{X}_t^y, \wt{Y}_{t-}^y, z)\mathcal{N}_{\mathbf{p}}(\d t,\d z), \quad \wt{X}_0^y=y,\, \wt{Y}_0^y=i_0.
  \end{cases}
  \end{equation}
  Note that under conditions (A1), (A2), SDEs \eqref{ac-1} and \eqref{ac-2}  both admit unique strong solution, which can be proved in the same way as in \cite[Theorem 2.3]{Sh15}.
  Then, according to Skorokhod's representation theorem (cf. \cite{YZ} or \cite[Theorem 2.2]{SWW}), $(X_t^x, Y_t^x)$ satisfies  \eqref{a-4.8} and $(\wt{X}_t^y, \wt Y_t^y)$ satisfies \eqref{a-4.9}.

 \noindent\textbf{Step 2}.  Based on the coupling process constructed above, we proceed to estimate the quantity
  $\dis \int_0^t\p( Y_s^x\neq \wt Y_s^y)\d s$ by induction.  It follows from the definition of $\mathcal{N}_{\mathbf{p}}(\d t, \d z)$, the estimate \eqref{U-1} that there exists a $\tilde c_1>0$  such that for $\delta>0$
  \[\p\big(\mathcal{N}_{\mathbf{p}}([0,\delta]\!\times\! U_{n })\geq 2\big)=1-\e^{-\mathbf{m}(U_n)\delta}-\mathbf{m}(U_n) \delta\e^{-\mathbf{m}(U_n)\delta}\!\leq \! \tilde c_1\delta^2,\quad n\geq 1.\]
  Then,
  \begin{align*}
    &\p(   Y_\delta^x\neq \wt Y_\delta^y|Y_0^x=\wt{Y}_0^y=i_0)\\
    &=\p\big(  Y_\delta^x\neq \wt Y_\delta^y, \mathcal{N}_{\mathbf{p}}([0,\delta] \!\times\! U_{i_0})=1|Y_0^x=\wt{Y}_0^y=i_0 \big)\\
    &\quad +\p\big(  Y_\delta^x\neq \wt Y_\delta^y,\mathcal{N}_{\mathbf{p}}([0,\delta]\!\times \! U_{i_0})\!\geq 2 |Y_0^x=\wt{Y}_0^y=i_0 \big) \\
    &\leq \p\big(  Y_\delta^x\neq \wt Y_\delta^y, \mathcal{N}_{\mathbf{p}}([0,\delta] \!\times\! U_{i_0})=1|Y_0^x=\wt{Y}_0^y=i_0 \big) +\tilde c_1\delta^2\\
    &=\int_0^\delta\!\p\big(\xi_1^{(i_0)}\!\in\! \!\bigcup_{j\in\S}\big(\Gamma_{i_0j}(X_s^x)\Delta \Gamma_{i_0j}(\wt{X}_s^y)\big),\tau_1^{(i_0)}\! \in\! \d s,\tau_2^{(i_0)}\!\geq \delta-s\big)+\tilde c_1 \delta^2,
  \end{align*} where $A\Delta B:=(A\backslash B )\cup (B\backslash A)$ for Borel sets $A$, $B$ in $\R$.
  Note that  for $s\leq \tau_1^{(i_0)}$, $X_s^x=X_s^{(i_0)}$ and $\wt{X}_s^y=\wt X_s^{(i_0)}$, where
  \begin{align*}
    X_s^{(i_0)}&=x+\int_0^sb(X_r^{(i_0)},i_0)\d r+\int_0^s\sigma(X_r^{(i_0)}, i_0)\d W_r,\\
    \wt X_s^{(i_0)}&=y+\int_0^sf(\wt{X}_r^{(i_0)}, i_0)\d r +\int_0^s g(\wt{X}_r^{(i_0)}, i_0)\d W_r.
  \end{align*}
  Therefore, due to the mutual independence of $\mathcal{N}_{\mathbf{p}}(\d t,\d z)$ and $(W(t))$, and the construction of $\Gamma_{ij}(z)$, we have
  \begin{equation}\label{ac-3}
  \begin{split}
    &\p(   Y_\delta^x\neq \wt Y_\delta^y|Y_0^x=\wt{Y}_0^y=i_0)\\
    &\leq \tilde c_1^2\delta^2+\int_0^\delta \E\Big[\sum_{j\neq i_0}|q_{i_0j}(X_s^x)-q_{i_0j}(\wt X_s^y)|\Big]\e^{-\mathbf{m}(U_{i_0})\delta}\d s\\
    &\leq \tilde c_1\delta^2 +\int_0^\delta\!\E\big[\|Q(X_s^x)-Q(\wt X_s^y)\|_{\ell_1}\big]\d s.
  \end{split}
  \end{equation}

  Now, let us consider $\p(Y_{2\delta}^x\neq \wt{Y}^y_{2\delta})$. It is clear that
  \begin{equation}\label{ac-4}
  \begin{split}
    &\p(Y_{2\delta}^x\neq \wt{Y}_{2\delta}^y)\\
    &=\p(Y_{2\delta}^x\neq \wt{Y}_{2\delta}^y|Y_\delta^x=\wt{Y}_\delta^y) \p(Y_{\delta}^x=\wt{Y}_\delta^y)+\p( Y_{2\delta}^x\neq \wt{Y}_{2\delta}^y, Y_\delta^x\neq \wt Y_\delta^y)\\
    &\leq \p(Y_{2\delta}^x\neq \wt{Y}_{2\delta}^y|Y_\delta^x=\wt{Y}_\delta^y)+ \p(Y_\delta^x\neq \wt{Y}_\delta^y).
  \end{split}
  \end{equation}
  Due to \eqref{ac-1} and \eqref{ac-2},
  \begin{equation}\label{ac-5}
  \begin{aligned}
    &\p(Y_{2\delta}^x\neq \wt Y_{2\delta}^y|Y_\delta^x= \wt Y_\delta^y)\\
    &\leq \p(Y_{2\delta}^x\neq \wt{Y}_{2\delta}^y,\mathcal{N}_{\mathbf{p}}((\delta, 2\delta])=1|Y_\delta^x=\wt Y_\delta^y) +\p( \mathcal{N}_{\mathbf{p}}((\delta, 2\delta])\geq 2|Y_\delta^x=\wt Y_\delta^y)\\
    &=\int_\delta^{2\delta}\!\!\p\big(\Delta p(s)\!\in\!\cup_{j\in\S}\big\{\Gamma_{Y_\delta^x j}(X_s^x)\Delta\Gamma_{\wt{Y}_\delta^y j}(\wt{X}_s^y)\big\},\tau_1^\delta\!\in\!\d s, \tau_2^\delta>2\delta-s\big)+\tilde c_1\delta^2,
  \end{aligned}
  \end{equation}
  where $\tau_1^\delta,\tau_2^\delta$ denote the first and second jump of $(\mathbf{p}(t))$ after time $\delta$. Note also that given $\F_{\delta}$, for $s\in[\delta, \tau_1^\delta]$, $X_s^x$ and $\wt {X}_s^y$ depend only on   $(W_r)_{r\in [\delta,s)}$. Based on the mutual independence of $(W_t)$ and $(\mathbf{p}(t))$, and their independent increment property, we   get from \eqref{ac-5} that
  \begin{equation}\label{ac-6}
  \begin{aligned}
     &\p(Y_{2\delta}^x\neq \wt Y_{2\delta}^y|Y_\delta^x=\wt Y_\delta^y) \leq \int_{\delta}^{2\delta}\!\!\E\big[\|Q(X_s^x)- Q(\wt X_s^y)\|_{\ell_1}\big] \d s+\tilde c_1\delta^2.
  \end{aligned}
  \end{equation}
  Inserting the estimates \eqref{ac-3} and \eqref{ac-6} into \eqref{ac-4}, we obtain that
  \begin{equation}\label{ac-7}
  \p(Y_{2\delta}^x\neq \wt Y_{2\delta}^y)\leq
  \int_0^{2\delta}\E\big[\|Q(X_s^x)-Q(\wt X_s^y)\|_{\ell_1}\big]\d s+2\tilde c_1\delta^2.
  \end{equation}
 Deducing inductively, we get
  \begin{equation}\label{ac-8}
  \p(  Y_{k\delta}^x\neq \wt Y_{k\delta}^y)\leq \int_0^{k\delta}\!\E\big[\|Q(X_s^x)-Q(\wt X_s^y)\|_{\ell_1}\big]\d s+k\tilde c_1\delta^2, \quad k\geq 3.
  \end{equation}

  Denote $N(t)=\big[\frac{t}{\delta}\big]$, the integer part of $t/\delta$, $t_k=k\delta$ for $k\leq N(t)$ and $t_{N(t)+1}=t$ for $t>0$. It follows from \eqref{ac-8} that
  \begin{align*}
    &\int_0^t\!\p(  Y_s^x\neq \wt Y_s^y)\d s\\
    &=\sum_{k=0}^{N(t)}\!\int_{t_k}^{t_{k+1}}\!\!\big\{ \p(  Y_s^x\neq \wt Y_s^y,   Y_{k\delta}^x=\wt Y_{k\delta}^y) \!+\! \p(  Y_s^x\neq \wt Y_s^y,  Y_{k\delta}^x\neq \wt Y_{k\delta}^y)\big\}\d s\\
    &\leq \sum_{k=0}^{N(t)}\!\int_{t_k}^{t_{k+1}}\!\!\p(  Y_s^x\neq \wt Y_s^y|  Y_{k\delta}^x=\wt Y_{k\delta}^y)\d s+\sum_{k=0}^{N(t)}\p(  Y_{k\delta}^x\neq \wt Y_{k\delta}^y) (t_{k+1}-t_k)\\
    &\leq  \sum_{k=0}^{N(t)}\!\int_{t_k}^{t_{k+1}}\!\! \p\big(\mathcal{N}_{\mathbf{p}} \big((t_k,t_{k+1}]\big)\!\geq\! 1\big )\d s+\!\delta\sum_{k=0}^{N(t)}\!\int_0^{k\delta}\!\! \Big( \E\big[ \|Q(X_s^x)\!-\!Q(\wt X_s^y)\|_{\ell_1}\big]\d s + k\tilde c_1\delta^2\Big)\\
    &\leq \big(1-\e^{-\kappa\delta}\big) t\!+\! \frac{(1+N(t))N(t)}{2}\tilde c_1\delta^3+\!\delta(N(t)+1)\int_0^t\!\E \big[\|Q(X_s^x)-Q(\wt X_s^y)\|_{\ell_1}\big]\d s.
  \end{align*}
  Letting $\delta\downarrow 0$, as $\delta (N(t)+1)\to t$, this yields that
  \[\frac1 t\int_0^t\!\p( Y_s^x\neq \wt Y_s^y)\d s \leq \int_0^t \!\E \big[\|Q(X_s^x)-Q(\wt X_s^y)\|_{\ell_1}\big]\d s.\]
  The proof of Lemma \ref{lem-1} is complete.
\fin

\section{Arguments of the main results}\label{argument}

This section is devoted to the arguments of the results presented in Section \ref{sec-2}.

We begin with proving the regularity of $x\mapsto \pi^x$ under strongly ergodic condition, which is based on the integration by parts formula for continuous-time Markov chains. The application of total variation norm and taking supremum in the initial value $i$ over $\S$ play an important role in the argument. 

\noindent\textbf{Argument of Proposition \ref{prop-1}} \quad
Using the integration by parts formula for continuous Markov chains (cf. \cite[Theorem 3.5]{Phi} or \cite[Theorem 13.40]{Chen}),
\begin{equation}\label{e-p-1}
P_t^yh(i)-P_t^xh(i)=\int_0^t\!P_{t-s}^y\big(Q(y)-Q(x)\big) P_s^x h(i) \d s,\quad t>0, \ h\in \mathscr{B}_b(\S).
\end{equation}
For any $h:\S\to \R$ with $|h|\leq 1$ and any $0\leq s\leq t$,
\begin{align*}
  \sup_{i\in\S} \big|P_{t-s}^y\big( Q(y)\!-\!Q(x)\big) P_s^x h(i)\big|&\leq \sup_{i\in \S} \big|(Q(y)-Q(x))P_s^x h(i)\big|\\
  &=\sup_{i\in \S}\big|(Q(y)-Q(x))P_s^x(h-\pi^x(h))(i)\big|\\
  &\leq 2\|Q(y)-Q(x)\|\sup_{i\in\S}|P_s^xh (i)-\pi^s(h)|\\
  &\leq 2\|Q(y)-Q(x)\| \sup_{i\in\S} \|P_s^x(i,\cdot)-\pi^x\|_\var,
\end{align*}
where, due to the conditions (A2) and (A3),  the operator norm
\begin{align*}
\|Q(y)-Q(x)\|&:=\sup\big\{|(Q(y)-Q(x))h(i)|; i\in \S, |h|\leq 1\big\}\\
&\leq  2\sup_{i\in\S} \sum_{j\in \S, j\neq i} |q_{ij}(y)-q_{ij}(x)|=2\|Q(y)-Q(x)\|_{\ell_1}\\
&\leq  2K_3|x-y|.
\end{align*}
Combining this estimate with (A4),  we get from \eqref{e-p-1}  that
\begin{equation}\label{e-p-2}
|P_t^yh(i)-P_t^xh(i)|\leq 4 c_1 K_3|x-y|\int_0^t\e^{-\lambda_1 s}\d s=\frac{4K_3 c_1}{\lambda_1} |x-y|\big(1-\e^{-\lambda_1 t}\big),
\end{equation}
and further
\begin{equation}\label{e-p-3}
\|P_t^y(i,\cdot)-P_t^x(i,\cdot)\|_{\var}\leq \frac{4K_3 c_1}{\lambda_1} |x-y|\big(1-\e^{-\lambda_1 t}\big)
\end{equation} by the arbitrariness of $h$ in \eqref{e-p-2}.

For any $h:\S\to \R$ with $|h|\leq 1$, it holds
\begin{equation}\label{e-p-4}
\begin{aligned}
\big|\pi^y(h)-\pi^x(h)\big|&=\Big|\sum_{i\in\S}  \pi_i^y P_t^y h(i)-\sum_{i\in \S} \pi_i^x P_t^x h(i)\Big|\\
  &\leq \sum_{i,j\in\S} \pi_j^y\pi_i^x\big|P_t^y h(j)-P_t^x h(i)\big|\\
  &\leq \sum_{j\in\S} \pi_j^y\big|P_t^y h(j)-P_t^x h(j)\big|+\sum_{i,j\in\S} \pi_j^y \pi_i^x \big| P_t^x h(j)-P_t^x h(i)\big|.
\end{aligned}
\end{equation}
By (A4), it holds
\begin{equation}\label{e-p-5}
|P_t^x h(i)-P_t^x h(j)|\leq |P_t^x h(j)-\pi^x(h)|+|P_t^x h(i)-\pi^x(h)|\leq 2c_1\e^{-\lambda_1 t}.
\end{equation}
Inserting \eqref{e-p-3}, \eqref{e-p-5} into \eqref{e-p-4}, we get
\begin{equation*}\label{e-p-6}
\begin{split}
  |\pi^y(h)-\pi^x(h)|&\leq \frac{4 K_3 c_1}{\lambda_1} |x-y|\big(1-\e^{-\lambda_1 t}\big)+2c_1\e^{-\lambda_1 t}.
\end{split}
\end{equation*}
Letting $t\to \infty$ and taking supremum over $h$ with $|h|\leq 1$, we obtain that
\[\|\pi^y-\pi^x\|_{\var}\leq \frac{4 K_3 c_1}{\lambda_1} |x-y|,\]
which is the desired conclusion, and the proof of Proposition \ref{prop-1} is completed.
\fin

As a direct application of Proposition \ref{prop-1}, it follows from the Lipschitz continuity of $b$ that $\bar{b}$ is also Lipschitz continuous. In fact,
\begin{equation}\label{a-6}
\begin{split}
  |\bar{b}(x)-\bar{b}(y)|&=\big|\sum_{i\in\S} b(x,i)\pi^x_i-\sum_{i\in\S}b(y,i)\pi_i^y\big|\\
  &\leq \sup_{i\in\S}|b(x,i)|\|\pi^x-\pi^y\|_\var +\big|\sum_{i\in\S}(b(x,i)-b(y,i))\pi_i^y\big|\\
  &\leq  \Big(K_1+\frac{2K_3 c_1}{\lambda_1}\Big) |x-y| , \qquad x,y\in\R^d.
\end{split}
\end{equation}

\noindent \textbf{Argument of Example \ref{exam-1}}\ Let $\mu_1^x=1$,
\[\mu_{n+1}^x=\frac{b_1b_2\ldots b_n}{a_2a_3\ldots a_{n+1}}=x^n ,\quad n\geq 1.\] Then,
\[\sum_{n=1}^\infty \mu_n^x=\frac{1}{1-x}<\infty,\quad \text{due to $x\in (0,1)$},\] and
\[ \sum_{n=1}^\infty \frac1{\mu_n^xb_n}\sum_{k=1}^n\mu_k^x= \sum_{n=1}^\infty \frac{1-x^n}{(1-x)x^n}=\infty. \]
According to the ergodic criterion for birth-death processes (cf. \cite[Chapter 1]{Chen}), the birth-death process $(Y_t^x)_{t\geq 0}$ is ergodic for every $x\in (0,1)$. Its invariant probability measure $\pi^x$ is given by $\pi_i^x=\frac{\mu_i^x}{\sum_{n=1}^\infty \mu_n^x}=(1-x)x^{i-1}$ for $i\geq 1$, which gives us \eqref{e-exam-2}.
Moreover,
one can check
\[\sup_{n\geq 2}\sum_{k=n}^\infty\mu_k^x \sum_{j\leq n-1} \frac1{\mu_j^xb_j}=  \sup_{n\geq 2}  \frac{1-x^{n-1}}{(1-x)^2} <\infty,\]
and hence $(Y_t^x)_{t\geq 0}$ is exponentially ergodic. However, by virtue of \cite[Theorem 3.1]{Mao02}, the birth-death process $(Y_t^x)_{t\geq 0}$ is not strongly ergodic since
\[\sum_{i=1}^\infty \frac{1}{\mu_i^xb_i}\sum_{j=i+1}^\infty \mu_j^x=\sum_{i=1}^\infty \frac{1}{1-x}=\infty.\] For the birth-death process $(Y_t^x)_{t\geq 0}$, its rate of exponential ergodicity is equivalent to the exponential $L^2$-convergence rate; see, \cite[Theorem 5.3]{Chen91}. Exponential $L^2$-convergence of Markov processes are closely related to the extensively studied Poincar\'e inequality and spectral gap of infinitesimal generators. There are many works devoted to the estimates of exponential $L^2$-convergence rate. Applying \cite[Example 5.7]{Chen91}, the exponential convergence rate of $(Y_t^x)_{t\geq 0}$ is given by
\begin{equation}\label{e-p-7}
\lambda(x)=\big(1-\sqrt{x} \big)^2,\qquad x\in (0,1).
\end{equation}

At last, we shall show that
\begin{equation}\label{e-p-8}
\sup_{x\neq y}\frac{\|\pi^x-\pi^y\|_\var}{|x-y|^\beta}=\infty,\quad \forall \beta\in (0,1],
\end{equation}
which yields \eqref{e-exam-3} and $x\mapsto \pi^x$ is not H\"older continuous with any exponent $\beta\in (0,1)$.

Indeed, we only need to consider the case $x>y$ in \eqref{e-p-8}. Due to the expression of $\pi^x$ in \eqref{e-exam-2},
consider the function $f(z)=(1-z) z^n$ on $(0,\infty)$. It holds
\[f'(z)=\Big(\frac n{n+1}-z\Big)(n+1)z^{n-1}.\] Therefore, when
$n>\frac{1-x}{x}>\frac{1-y}{y}$, $f'(z)>0$ for all $z\in [y,x]$. This implies that $\pi_{n+1}^x=(1-x)x^n>(1-y)y^n=\pi_{n+1}^y$.
Let $n_x=\inf\{m\in \N; m\geq (1-x)/x\}$.
Therefore,
\begin{equation}\label{e-p-9}\|\pi^x-\pi^y\|_\var =\sum_{n=1}^\infty |\pi_n^x-\pi_n^y|\geq \sum_{n=n_x+1}^\infty\big( \pi_n^x-\pi_n^y\big)= x^{n_x}-y^{n_x}.
\end{equation}

Take $x=1-\frac 1m$ and $y=\frac{2m-2}{2m-1} x$ for $m\geq 2$, then $1>x>y>0$ and $n_x=m$.
For any $\beta\in (0,1]$, due to \eqref{e-p-9},
\begin{equation}\label{e-p-10}
\begin{split}
  \sup_{x\neq y}\frac{\|\pi^x-\pi^y\|_\var}{|x-y|^\beta}&\geq \lim_{m\to \infty} \frac{(1-\frac 1m)^m-(1-\frac1m)^m (1-\frac1{2m-1})^m}{(1-\frac 1m)^\beta \frac 1{(2m-1)^\beta}}\\
  &=\lim_{m\to\infty} (2m-1)^\beta \Big(1-(1-\frac1{2m-1})^m\Big)\\
  &=\infty.
\end{split}
\end{equation}
All assertions in Example \ref{exam-1} have been proved.
\fin

Before presenting the proofs of Theorem \ref{thm-1} and Theorem \ref{thm-1.5}, let us introduce the main challenge in the proofs.
Firstly,
we should pay more attention to the difficulty caused by the full dependence of the two time-scale system $(X_t^\ea,Y_t^\ea)$.  To overcome this difficulty,  we shall use the coupling method developed in Section \ref{coupling}. Secondly,
we need to pay attention to the essential difference  between the distributions of $(X_t^\ea)_{t\in [0,T]}$ and those of $(Y_t^\ea)_{t\in [0,T]}$  for $\veps,\alpha\in (0,1)$ given $T>0$. Precisely, for each fixed $T>0$, let $\mathcal{C}([0,T];\R^d)$ be the space of continuous paths from $[0,T]$ to $\R^d$, and $\mathcal{D}([0,T];\S)$ the Skorokhod space containing right continuous paths with left limits. Then under condition (A1),   the distributions of $\{(X_t^\ea)_{t\in [0,T]};\veps,\alpha>0\} $ in $\mathcal{C}([0,T];\R^d)$ is tight. However, the distributions of  $\{(Y_t^\ea)_{t\in [0,T]};\veps,\alpha>0\}$ in $\mathcal{D}([0,T];\S)$ is not tight, which can be seen from the following simple and meaningful example given in \cite[Example 7.3, p.172]{YZ98}.
\begin{myexam}[\cite{YZ98}]\label{exam-2}
Let $(\La_t^\alpha)_{t\in [0,T]}$ be a continuous time Markov chain on the state space $\S=\{1,2\}$ with transition rate
\[\frac1\alpha \begin{pmatrix}
  -\lambda&\lambda\\ \mu &-\mu
\end{pmatrix},\]
for some $\lambda,\mu>0$. Then, for each $T>0$ the collection of distributions of $(\La_t^\alpha)_{t\in [0,T]}$ for $\alpha\in (0,1)$ is not tight.
\end{myexam}

\noindent \textbf{Argument of Theorem \ref{thm-1}}
Let $(X_t^\ea, Y_t^\ea)$ be a solution to SDEs  \eqref{o-1}, \eqref{o-2}. Based on Skorokhod's representation theorem, similar to SDE  \eqref{ac-1}, $(X_t^\ea,Y_t^\ea)$ can be expressed as a solution to SDEs driven by a Brownian motion and a Poisson random measure respectively. In the following, this expression of $(X_t^\ea, Y_t^\ea)$ helps us to use the method introduced in Section \ref{coupling} to construct the desired coupling process so as to decouple the interaction between $(Y_t^\ea)$ and $(X_t^\ea)$.

For $\delta>0$, let $t(\delta)=[\frac t\delta]\delta$, where $[\frac t\delta]=\max\{n\in \N; n\leq \frac t \delta\}$.   Due to the boundedness of $b$ and $\sigma$ in (A1), it follows from  \eqref{o-1} that
\begin{align*}
  \E|X_t^\ea-X_{t(\delta)}^\ea|&\leq \E\int_{t(\delta)}^t\!|b(X_s^\ea,Y_s^\ea)|\d s +\!\sqrt{\veps}\Big(\E \int_{t(\delta)}^t\! \|\sigma(X_s^\ea,Y_s^\ea)\|^2\d s\Big)^{\frac 12} \\
  &\leq K_2(\delta+\sqrt{\veps \delta}).
\end{align*}
Using the triangle inequality, we divide the estimate of $\E|X_t^\ea-\bar{X}_t|$ into five terms:
  \begin{equation}\label{a-7}
  \begin{split}
    &\E|X_t^\ea-\bar{X}_t|\\
    &\leq \E \Big|\int_0^t\!b(X_s^\ea,Y_s^\ea)-\bar b(\bar{X}_s)\d s\Big|+\sqrt{\veps}\E\Big|\int_0^t\! \sigma(X_s^\ea,Y_s^\ea)\d W_s\Big|\\
    &\leq \E\Big|\int_{t(\delta)}^t\!b(X_s^\ea,Y_s^\ea)-\bar b (\bar{X}_s)\d s\Big|+\E\Big|\int_0^{t(\delta)}\!b(X_s^\ea,Y_s^\ea) -b(X_{s(\delta)}^\veps,Y_s^\ea)\d s\Big|\\
    &\quad +\E \Big| \int_0^{t(\delta)}\! b(X_{s(\delta)}^\veps,Y_s^\ea) -\bar{b}( X^\veps_{s(\delta)})\d s\Big|+\E\Big|\int_0^{t(\delta)} \bar{b}(X_{s(\delta)}^\veps)-\bar{b}(\bar{X}_s)\d s\Big|\\
    &\quad +\sqrt{\veps} \E\Big|\int_0^t\sigma(X_s^\ea, Y_s^\ea)\d W_s\Big|=:(I)+(I\!I)+(I\!I\!I)+(I\!V)+(V).
  \end{split}
  \end{equation}
   We shall estimate the right hand side of \eqref{a-7} terms by terms. By (A1)
   \begin{equation}\label{a-7.5}
   \begin{aligned}
   (I)=\E\Big|\int_{t(\delta)}^t b(X_s^\ea,Y_s^\ea)-\bar b(\bar{X}_s)\d s\Big|&\leq   2 K_2\delta,\\
   (I\!I)=\!\E\Big|\int_0^{t(\delta)}\!\! \! b(X_s^\ea,Y_s^\ea)\!-\!b(X_{s(\delta)}^\ea, Y_s^\ea)\d s\Big|& \leq K_1\!\!\int_0^{t(\delta)}\!\! \E|X_s^\ea\!-\!X_{s(\delta)}^\ea|\d s\\
   &\leq K_1 K_2 t\big( \delta\!+\!   \sqrt{\veps\delta} \big),
   \end{aligned}
   \end{equation}
   and \[ (V)=\sqrt{\veps}\E\Big|\int_0^t\sigma(X_s^\ea,Y_s^\ea)\d W_s\Big|\leq \sqrt{\veps} \Big(\E\int_0^t\|\sigma(X_s^\ea,Y_s^\ea)\|^2\d s\Big)^{1/2}\leq K_2\sqrt{\veps t}.\]

   To deal with term $(I\!I\!I)$, we divide the integral over $[0,t(\delta))$ into the integrals over subintervals $[k\delta, (k+1)\delta)$ via the following inequality
   \[\E\Big|\int_0^{t(\delta)}\!\! b(X_{s(\delta)}^\ea,Y_s^\ea)- \bar{b}(X_{s(\delta)}^\ea)\d s\Big|\leq \sum_{k=0}^{t(\delta)/\delta-1}\!\int_{k\delta}^{(k+1)\delta} \E \big|b(X_{s(\delta)}^\ea,Y_s^\ea)- \bar{b}(X_{k\delta}^\ea)\big|\d s,\]
   then at each subinterval $[k\delta, (k+1)\delta)$, $k\in \Z_+$ with $0\leq k\leq t(\delta)/\delta-1$, we introduce an auxiliary process $( \wt{Y}_t^{(k)})_{t\geq k\delta}$ constructed as in Lemma \ref{lem-1} such that:
   \begin{itemize}
   \item[$\mathrm{(i)}$]  Under the conditional expectation $\E\big[\,\cdot\,\big|\F_{k\delta}\big]$, $(\wt Y_t^{(k)})_{t\geq k\delta}$ is a Markov chain in $\S$ with transition rate matrix $(\frac 1{\alpha}q_{ij}(X_{k\delta}^\ea))_{i,j\in\S}$ and satisfies $\wt Y_{k\delta}^{(k)}=Y_{k\delta}^\ea$.
   \item[$\mathrm{(ii)}$] The following estimate holds: for $t>k\delta$,
       \begin{equation}\label{a-7.6}
       \frac{1}{t-k\delta}\int_{k\delta}^t\! \E\big[\mathbf1_{\{Y_s^\ea\neq \wt{Y}_s^{(k)}\}}\big|\F_{k\delta}\big]\d s\leq  \int_{k\delta}^t \E\big[\|Q(X_s^\ea)-Q(X_{k\delta}^\ea )\|_{\ell_1}\big|\F_{k\delta}\big] \d s.
       \end{equation}
   \end{itemize}
   Noting the   scaling $1/\alpha$ in the transition rate matrix of $(\wt{Y}_t^{(k)})_{t\geq k\delta}$, we have
   \[\E\big[f(\wt{Y}_{k\delta+s}^{(k)})\big|\F_{k \delta}\big]=P_{s/\alpha}^{X_{k\delta}^\ea}(f) (Y_{k\delta}^{\ea}) \] for any bounded function $f$ on $\S$, where $P_t^x$ denotes the semigroup corresponding to the $Q$-matrix $(q_{ij}(x))_{i,j\in\S}$ as before.
     By (A4), for any $h\in \B(\S)$ with $|h| \leq 1$, $ s> k\delta$,
   \begin{equation}\label{a-8}
   \begin{split}
    \E\big[\big| h(\wt Y_s^{(k)})- \pi^{X_{k\delta}^\veps}(h)\big|\big| \F_{k\delta}\big]&\leq \sup_{x\in\R^d}\sup_{i\in\S} \|P_{\frac{s-k\delta}{\alpha} }^{x}(i,\cdot)-\pi^{x} \|_{\var} \leq c_1\e^{-\lambda_1 \frac{s-k\delta}{\alpha}}.
    \end{split}
    \end{equation}
   Hence,
   \begin{equation*}
   \begin{split}
     &\int_{k\delta}^{(k+1)\delta}\!\!\E\big|b(X_{k\delta}^\ea, Y_s^\ea)-\bar{b}(X_{k\delta}^\ea)\big|\d s\\
     &\leq\!  \int_{k\delta}^{(k+1)\delta} \!\!\!\E\Big\{\E\Big[\big|b(X_{k\delta}^\ea ,Y_s^\ea)\!-\!b(X_{k\delta}^\ea,\wt Y_{s}^{(k)})\big| \!+\! \big| b(X_{k\delta}^\ea,\wt Y_{s}^{(k)}) \!-\!\bar{b}(X_{k\delta}^\ea)\big|\Big| \F_{k\delta}\Big]\Big\}\d s\\
     &\leq 2K_2 \alpha \!\int_{0}^{\frac{\delta}{\alpha}} \! \p(Y_{k\delta+\alpha r}^\ea \neq \wt{Y}_{k\delta+\alpha r}^{(k)})\d r \\
     &\quad +\!\alpha\!\int_{0}^{\frac{\delta}\alpha}\! \E\Big\{\E \Big[|b(X_{k\delta}^\ea, \wt{Y}_{k\delta+\alpha r}^{(k)})\!-\!\sum_{i\in \S}b(X_{k\delta}^\ea,i) \pi_i^{X_{k\delta}^\ea} \big|\Big|\F_{k\delta} \Big] \Big\}\d r
     \end{split}
     \end{equation*}
     \begin{equation}\label{a-9}
     \begin{split}
     &\leq  2K_2K_3\delta \int_0^{\frac{\delta}{\alpha}} \E\big[|X_{k\delta+\alpha r}^\ea -X_{k\delta}^\ea|\big]\d r + \alpha K_2 \int_0^{\frac\delta\alpha} \sup_{x\in\R^d}\sup_{i\in\S}\|P_r^x(i,\cdot) -\pi^x\|_{\var} \d r\\
     &\leq   2K_2^2K_3\frac{\delta}{\alpha}\int_0^\delta (r +\sqrt{\veps r})\d r+\frac{\alpha K_2c_1}{\lambda_1}\big(1-\e^{-\lambda_1 \delta/\alpha}\big)\\
     &\leq  2K_2^2K_3 \frac{\delta}{\alpha} \Big(\frac{\delta^2}{2}+ \frac{2\veps^{1/2}\delta^{3/2}}{3}\Big) +\frac{\alpha K_2c_1}{\lambda_1}\big(1-\e^{-\lambda_1 \delta/\alpha}\big),
   \end{split}
   \end{equation}
   where in the second inequality we used \eqref{a-7.6} and (A3), and in the third inequality we used \eqref{a-8}. Therefore,
   \begin{equation}\label{o-10}
   \begin{split}
  (I\!I\!I)= &\E\Big|\int_0^{t(\delta)}\!b(X_{s(\delta)}^\ea, Y_s^\ea)-\bar{b}(X_{s(\delta)}^\ea)\d s \Big| \\
   &\leq   2K_2^2K_3 \frac{t(\delta)}{\delta} \frac{\delta}{\alpha}\big(\frac{\delta^2}{2} +\frac{2 \veps^{1/2}\delta^{3/2}}{3}\big)\!+\!  t(\delta) \frac{\alpha}{\delta}\frac{  K_2c_1}{\lambda_1}\big(1-\e^{-\lambda_1 \delta/\alpha}\big)\\
   &\leq  2K_2^2K_3 t\frac{\delta}{\alpha}\big(\frac\delta 2+\frac{2\veps^{1/2} \delta^{1/2}}{3}\big)+  t\frac{\alpha}{\delta}\frac{  K_2c_1}{\lambda_1}\big(1-\e^{-\lambda_1 \delta/\alpha}\big).
   \end{split}
   \end{equation}
   Taking $\delta=\alpha^{3/4}$ and invoking \eqref{a-7.5}, we obtain that
   \begin{equation}\label{o-11}
   \begin{split}
     &\E\Big|\int_0^{t }\!b(X_{s(\delta)}^\ea, Y_s^\ea)-\bar{b}(X_{s(\delta)}^\ea)\d s \Big|\\ &
     \leq 2K_2\alpha^{3/4}\!+\! 2K_2^2K_3 t \alpha^{\frac 18}\big(\frac{\alpha^{\frac 38}}{2}+ \frac{2\veps^{\frac 12}}3\big)+ t\alpha^{\frac14} \frac{  K_2c_1}{\lambda_1}\big(1-\e^{-\lambda_1/ \alpha^{\frac 14}}\big).
   \end{split}
   \end{equation}
By \eqref{a-6},
\begin{equation}\label{o-12}
\begin{split}
(I\!V)&=\E\int_0^{t(\delta)}\!\!|\bar{b}(X_{s(\delta)}^\ea) -\bar{b} (\bar{X}_s)|\d s\leq (K_1+\frac{2K_3c_1}{\lambda_1})\int_0^{t(\delta)}\!\! \E|X_{s(\delta)}^\ea-\bar X_s|\d s\\
&\leq (K_1+\frac{2K_3c_1}{\lambda_1})K_2\big(\delta+\sqrt{\veps \delta}\big) t+(K_1+\frac{2K_3c_1}{\lambda_1})\int_0^t\E|X_s^\ea-\bar{X}_s|\d s.
\end{split}
\end{equation}
Consequently, inserting above estimates \eqref{o-12}, \eqref{o-11} into \eqref{a-7} by taking $\delta=\alpha^{3/4}$, we obtain that
\begin{equation}\label{o-13}
\E|X_t^\ea-\bar{X}_t|\leq \phi(\veps,\alpha)+(K_1+\frac{2K_3c_1}{\lambda_1} )\int_0^t\E|X_s^\ea-\bar X_s|\d s,
\end{equation}
where
\begin{align*}
\phi(\veps,\alpha)&= 2K_2\alpha^{\frac34}+K_1K_2 t(\alpha^{\frac34}+\veps^{\frac 12}\alpha^{\frac 38})+K_2\veps^{\frac 12} t^{\frac 12}+ 2K_2^2K_3 t\alpha^{\frac 18}\big(\frac{\alpha^{\frac38}}{2} +\frac 23\veps^{\frac 12}\big)\\
&\quad +K_2t \alpha^{\frac 14}\int_0^\infty \eta_r\d r +
(K_1+\frac{2K_3c_1}{\lambda_1} )K_2 t(\alpha^{\frac 34}+\veps^{\frac 12}\alpha^{\frac 38}).
\end{align*}
Applying Gronwall's inequality, we finally get
\[\lim_{\veps,\alpha\to 0}
\E|X_t^\ea-\bar{X}_t|\leq \lim_{\veps,\alpha\to 0} \phi(\veps,\alpha)\e^{(K_1+\frac{2K_3c_1}{\lambda_1})t}=0,\]
and the proof of this theorem is complete.
\fin

\noindent\textbf{Argument of Proposition \ref{prop-2}}\ \
For any bounded function $h$ on $\S$ with $|h|\leq 1$, take some $i_0\in\S$, and then it holds that
  \begin{align*}
    &|\pi^x(h)-\pi^y(h)|=\big|\sum_{i\in\S} h_i\pi_i^x\!-\!\sum_{i\in\S} h_i\pi_i^y \big|\\
    &\leq \big|\pi^x(h)\!-\!\frac 1t\int_0^tP_s^xh(i_0)\d s\big|+\big|\pi^y(h)\!-\!\frac 1t\int_0^t\!P_s^y h(i_0)\d s\big|+\big|\frac 1t\int_0^t\! \!\big(P_s^xh(i_0)-P_s^yh(i_0)\big)\d s\big|\\
    &\leq \frac 1t \int_0^t\!\big|P_s^x h(i_0)-\pi^x(h)\big|\d s+\frac 1t\int_0^t\!\big|P_s^y h(i_0)-\pi^y(h)\big|\d s+\frac 1t\int_0^t\!2\p(\wt Y_s^x\neq \wt Y_s^y)\d s\\
    &\leq \frac{1}t\!\int_0^t \!\! \|P_s^x(i_0,\cdot)-\pi^x\|_\var\d s+\frac{1}t\!\int_0^t\!\!  \|P_s^y(i_0,\cdot)-\pi^y\|_\var \d s\!+\!  \|Q(x)-Q(y)\|_{\ell_1}t\\
    &\leq \frac{2\theta(i_0)} t\int_0^\infty\eta_s\d s+ \|Q(x)-Q(y)\|_{\ell_1}t,\qquad \forall\,t>0,
  \end{align*}
  where we have used Lemma \ref{lem-1} and (A4), which ensures that $\int_0^\infty\eta_s\d s<\infty$. Then, by taking $t=\sqrt{\frac{2\theta(i_0)\int_0^\infty \eta_s\d s}{\|Q(x)-Q(y)\|_{\ell_1}}}$, we arrive at
  \begin{equation}\label{aa-6}
  |\pi^x(h)-\pi^y(h)|\leq \sqrt{2\theta(i_0)\int_0^\infty\!\! \eta_s\d s\,\|Q(x)-Q(y)\|_{\ell_1}}.
  \end{equation}
  By the arbitrariness of $h$ and   (A3),
  \[\|\pi^x -\pi^y\|_{\var}\leq  \Big(2 K_3\theta(i_0)\int_0^\infty \!\!\eta_s\d s \Big)^{\frac 12} \sqrt{|x-y|},\]
  and further the desired estimate \eqref{aa-5} by taking the infimum for $\theta(i_0)$ over $i_0\in\S$.\fin

Analogous to the deduction of \eqref{a-6}, under conditions (A1)-(A3)  and (A5), $\bar{b}$ is $1/2$-H\"older continuous by virtue of Proposition \ref{prop-2}.
According to Peano's theorem, ODE \eqref{a-4} must admit  a solution. However, it may loss the uniqueness of solution. Moreover, in contrast to the $L^1$-convergence in Theorem \ref{thm-1} in the strongly ergodic condition, we can only prove the weak convergence of $(X_t^\ea)$ to $(\bar{X}_t)$.

\noindent\textbf{Proof of Theorem \ref{thm-1.5}}  \ \ Denote by $\mathscr{L}^\ea$ the generator of $(X_t^\ea)$ given by
\begin{equation}\label{ab-0}
 \mathscr{L}^\ea f(x,i)\!=\!\la \nabla f(x), b(x,i)\raa \!+\!\frac \veps 2\mathrm{tr}\big((\sigma\sigma^\ast)(x,i)\nabla^2 f(x)\big),\  f\! \in\! C^2_b(\R^d), \, x\!\in \! \R^d,\,i\!\in\!\S.
\end{equation}
Here, for a matrix $A$, $A^\ast $ denotes its transpose and $\mathrm{tr}(A)$ its trace.
Let $T>0$ be fixed. Let $\mathcal{C}([0,T];\R^d)$ be endowed with uniform norm, i.e. $\|x_\cdot-y_\cdot\|_{\infty}=\sup_{t\in [0,T]}|x_t-y_t|$ for $x_\cdot, \, y_\cdot\in C([0,T];\R^d)$. Denote by $\mathcal{L}_{X^\ea}$ the law of the process $(X_t^\ea)_{t\in [0,T]}$ in the path space $\mathcal{C}([0,T];\R^d)$.

Due to the boundedness of $b$ and $\sigma$ in (A1), it is standard to show
\begin{equation}\label{ab-1}
\E\big[\sup_{t\in[0,T]}|X_t^\ea|^p\big]\leq C(T,x_0,p),\quad \forall \,p\geq 1,
\end{equation}
where $x_0=X_0^\ea$, $C(T,x_0,p)$ is a constant depending on $T,\,x_0$ and $p$.
By It\^o's formula, for $0\leq s<t\leq T$,
\begin{align*}
  \E|X_t^\ea-X_s^\ea|^4&\leq 8\E\Big|\int_s^t\!b(X_r^\ea,Y_r^\ea)\d r\Big|^4+\!8\veps^2\E\Big|\int_s^t\!\sigma( X_r^\ea,Y_r^\ea)\d W_r\Big|^4\\
  &\leq 8(t-s)^3\E\int_s^t\!|b(X_r^\ea, Y_r^\ea)|^4\d r+288 \veps^2(t-s)\E\int_s^t\! |\sigma(X_r^\ea, Y_r^\ea)|^4\d r\\
  &\leq C(t-s)^2
\end{align*} for some constant $C>0$. Combing this with $X_0^\ea=x_0$, the collection of laws  $ \mathcal{L}_{X^\ea}$ for $\veps,\,\alpha>0 $ over the space $\mathcal{C}([0,T];\R^d)$ is tight by virtue of \cite[Theorem 12.3]{Bill}. As a consequence, there is a subsequence $\{\mathcal{L}_{X^{\veps',\alpha'}}; \veps',\alpha'>0\}$ and a limit law $\mathcal{L}_{\wt X}$ in $\mathcal{C}([0,T];\R^d)$ such that $\mathcal{L}_{X^{\veps',\alpha'}}$ converges weakly to $\mathcal{L}_{\wt X}$ as $\veps',\,\alpha'\to 0$. According to Skorokhod's representation theorem with a slight abuse of notation, we may assume that $(X_t^{\veps',\alpha'})_{t\in [0,T]}$ converges almost surely to some $(\wt X_t)_{t\in[0,T]}$  in $\mathcal{C}([0,T];\R^d)$ as $\veps',\alpha'\to 0$.

In order to characterize the limit, we shall show that for any $f\in C_c^2(\R^d)$, the space of functions with compact support and continuous second order derivatives.
\begin{gather*}
  \text{$f(\wt X_t)-f(x_0)-\int_0^t \mathscr{L} f(\wt X_s)\d s$ is a martingale},
\end{gather*}
where
\begin{equation}\label{ab-2}
\mathscr{L}f(x)=\la\nabla f(x), \bar{b}(x)\raa,\quad \text{and}\ \ \bar{b}(x)=\sum_{i\in\S} b(x,i)\pi_i^x.
\end{equation}
This means that $(\wt X_t)$ is a solution to ODE \eqref{a-4}.

To this end, it suffices to show that for any $0\leq s<t\le T$, for any bounded $\F_s$ measurable function $\Phi$,
\begin{equation}\label{ab-3}
\E\Big[\Big(f(\wt X_t)-f(\wt X_s)-\int_s^t\!\mathscr{L}f(\wt X_r)\d r\Big)\Phi\Big]=0,\quad \forall\, f\in C_c^2(\R^d).
\end{equation}
As a solution to SDE \eqref{o-1}, $(X_t^{\veps',\alpha'})$ satisfies
\begin{equation}\label{ab-4}
\E\Big[\Big( f(X_t^{\veps',\alpha'})- f(X_s^{\veps',\alpha'}) -\int_s^t\! \mathscr{L}^{\veps',\alpha'}f(X_r^{\veps', \alpha'}, Y_r^{\veps',\alpha'})\d r\Big)\Phi\Big]=0.
\end{equation}
By the dominated convergence theorem, it is clear that
\[\lim_{\veps',\alpha'\to 0} \E\Big[\big(f(X_t^{\veps',\alpha'})- f(X_s^{\veps',\alpha'})\big)\Phi\Big]= \E\Big[\big(f(\wt X_t)-f(\wt X_s)\big)\Phi\Big].
\]
Hence, to derive \eqref{ab-3} from \eqref{ab-4} we only need to show
\[\lim_{\veps',\alpha'\to 0} \E\Big[\int_s^t\! \big(\mathscr{L}^{\veps',\alpha'} f(X_r^{\veps',\alpha'}, Y_r^{\veps',\alpha'})-\mathscr{L}f(\wt X_r)\big)\d r\Big|\F_s\Big]=0.\]
According to the expression  \eqref{ab-0}, \eqref{ab-2} of $\mathscr{L}^{\veps',\alpha'}$, $\mathscr{L}$ and the boundedness of $\sigma$, it suffices to show
\begin{equation}\label{ab-5}
\lim_{\veps',\alpha'\to 0}
\E\Big[\int_s^t\!\big(\la \nabla f(X_r^{\veps',\alpha'}), b(X_r^{\veps',\alpha'}, Y_r^{\veps',\alpha'})\raa-\la \nabla f(\wt X_r), \bar{b}(\wt X_r)\raa \big)\d r \Big|\F_s\Big]=0.
\end{equation}
Similar to the treatment of \eqref{a-7}, we shall use the time discretization method and the coupling method to show \eqref{ab-5}.

Precisely, for $\delta>0$, let $r(\delta)=s+\big[\frac{r-s}{\delta}\big] \delta $ for $r\in [s, t]$.
\begin{equation}
\begin{split}
&\E\Big[ \int_s^t\!\big| \la \nabla f(X_r^{\veps',\alpha'}), b(X_r^{\veps',\alpha'}, Y_r^{\veps',\alpha'})\raa -\la \nabla f(\wt X_r), \bar{b}(\wt X_r)\raa \big|\d r \Big|\F_s\Big]\\
&\leq \E\Big[\int_s^t\!\Big(\big|\la \nabla f(X_r^{\veps',\alpha'}), \nabla b(X_r^{\veps',\alpha'}, Y_r^{\veps',\alpha'})\raa -\la\nabla f(X_{r(\delta)}^{\veps',\alpha'}),
b(X_{r(\delta)}^{\veps',\alpha'},  Y_r^{\veps',\alpha'})\raa \big|\\
&\qquad\quad +\big| \la \nabla f(X_{r(\delta)}^{\veps',\alpha'}), b(X_{r(\delta)}^{\veps',\alpha'}, Y_r^{\veps',\alpha'})\!-\!\bar{b} (X_{r(\delta)}^{\veps',\alpha'})\raa \big| \\ &\qquad\quad  +\big|\la \nabla f(X_{r(\delta)}^{\veps',\alpha'}), \bar{b}(X_{r(\delta)}^{\veps',\alpha'})\raa-\la \nabla f(\wt X_r),\bar{b}(\wt X_r)\raa \big|\Big)\d r\Big|\F_s\Big]\\
&=:  \Upsilon_1 +\Upsilon_2+\Upsilon_3.
\end{split}
\end{equation}
Applying the boundedness and Lipschitz continuity of $b$, the H\"older continuity of $\bar{b}$,  and the almost sure convergence of $(X_t^{\veps',\alpha'})_{t\in [0,T]}$ to $(\wt X_t)_{t\in [0,T]}$ as $\veps',\alpha'\to 0$, we obtain that
\begin{align*}
  &\lim_{\veps',\alpha',\delta\to 0} \big(\Upsilon_1+\Upsilon_3\big)\\
  &=\lim_{\veps',\alpha',\delta\to 0}  \E\Big[\int_s^t\!\Big(\big|\la \nabla f(X_r^{\veps',\alpha'}), \nabla b(X_r^{\veps',\alpha'}, Y_r^{\veps',\alpha'})\raa -\la\nabla f(X_{r(\delta)}^{\veps',\alpha'}),
b(X_{r(\delta)}^{\veps',\alpha'},  Y_r^{\veps',\alpha'})\raa \big|\\
&\qquad \quad \qquad \quad+  \big|\la \nabla f(X_{r(\delta)}^{\veps',\alpha'}), \bar{b}(X_{r(\delta)}^{\veps',\alpha'})\raa-\la \nabla f(\wt X_r),\bar{b}(\wt X_r)\raa \big|\Big)\d r\Big|\F_s\Big]=0.
\end{align*}
Meanwhile, let $N_t =[(t-s)/\delta]$, $s_k=s+k\delta$ for $0\leq k\leq N_t$ and $s_{N_t+1}=t$.
It holds
\begin{equation}\label{ab-6.5}
  \Upsilon_2 \leq \sum_{k=0}^{N_t}\E\Big[\int_{s_k}^{s_{k+1}}\! \big|\la \nabla f(X_{s_k}^{\veps',\alpha'}), b(X_{s_k}^{\veps',\alpha'}, Y_r^{\veps',\alpha'} ) -\bar{b}(X_{s_k}^{\veps',\alpha'}) \raa \big|\d r \Big|\F_s\Big].
\end{equation}

Next, for each $0\leq k\leq N_t$, let us construct an auxiliary Markov chain $(\wt{Y}_t^{(k)})_{t\geq s_k}$ in the way as Lemma \ref{lem-1} such that under the conditional expectation $\E[\,\cdot\,|\F_{s_k}]$, $(\wt{Y}_t^{(k)})_{t\geq s_k}$ is a Markov chain in $\S$ with transition rate $\big(\frac1{\alpha'}q_{ij} (X_{s_k}^{\veps',\alpha'})\big)_{i,j\in\S}$ and satisfies $\wt{Y}_{s_k}^{(k)}=Y_{s_k}^{\veps',\alpha'}$.
Moreover, it also satisfies
\begin{equation}\label{ab-7}
\frac{1}{t-s_k}\int_{s_k}^t\! \E\big[\mathbf1_{\{Y_r^{\veps',\alpha'}\neq \wt{Y}_r^{(k)}\}}\big|\F_{s_k}\big]\d r\leq \int_{s_k}^t \E\big[\|Q(X_r^{\veps',\alpha'})-Q(X_{s_k}^{\veps',\alpha'} )\|_{\ell_1}\big|\F_{s_k}\big] \d r.
\end{equation}
Then,
\begin{align*}
  &\E\Big[\int_{s_k}^{s_{k+1}}\! \big|\la \nabla f(X_{s_k}^{\veps',\alpha'}), b(X_{s_k}^{\veps',\alpha'}, Y_r^{\veps',\alpha'})- \bar{b}(X_{s_k}^{\veps',\alpha'}) \raa\big|\d r\Big|\F_{s_k}\Big]\\
  &\leq \int_{s_k}^{s_{k+1}}\!\!\!   \|\nabla f \|_\infty \Big(\E\big[|b(X_{s_k}^{\veps',\alpha'}, Y_r^{\veps',\alpha'})-b(X_{s_k}^{\veps',\alpha'}, \wt{Y}_r^{(k)})|\big|\F_{s_k}\big]\\ &\qquad\quad  +\!\E\big[|b(X_{s_k}^{\veps',\alpha'},
   \wt{Y}_r^{(k)})-\bar{b}(X_{s_k}^{\veps',\alpha'}) | \big|\F_{s_k}\big]\Big)\d r\\
  &\leq 2K_2\|\nabla f\|_{\infty}\alpha' \int_0^{\frac{\delta}{\alpha'}}\!\E\Big[\mathbf1_{ Y_{s_k+\alpha' r}^{\veps',\alpha'}\neq \wt{Y}_{s_k+\alpha' r}^{(k)}}\big|\F_{s_k}\Big]\d r \\ &\quad +\!K_2\|\nabla f\|_\infty \alpha'\int_0^{\frac{\delta}{\alpha'}}\!\!\E \Big[ \|P_{r}^{X_{s_k}^{\veps',\alpha'}}(Y_{s_k}^{\veps', \alpha'}, \cdot)\!-\! \pi^{X_{s_k}^{\veps',\alpha'}} \|_{\var}\big|\F_{s_k}\Big] \d r,
\end{align*} where $\|\nabla f\|_\infty=\!\sup_{z\in\R^d}\! |\nabla f(z)|<\infty$ for $f\in\! C_c^2(\R^d)$.
By virtue of \eqref{ab-7}, (A3) and (A5),
\begin{equation}\label{ab-8}
\begin{aligned}
 &\E\Big[\int_{s_k}^{s_{k+1}}\! \big|\la \nabla f(X_{s_k}^{\veps',\alpha'}), b(X_{s_k}^{\veps',\alpha'}, Y_r^{\veps',\alpha'})- \bar{b}(X_{s_k}^{\veps',\alpha'}) \raa\big|\d r\Big|\F_{s_k}\Big]\\
 &\leq 2K_2K_3\|\nabla f\|_\infty \delta\!\int_0^{\frac{\delta}{\alpha'}}\!\! \E\big[|X_{s_k+\alpha'r}^{\veps',\alpha'} \!-\! X_{s_k}^{\veps', \alpha'}|\big|\F_{s_k}\big]\d r \!+\! \alpha'K_2\|\nabla f\|_\infty\theta(Y_{s_k}^{ \veps',\alpha'}\!)\!  \int_0^{ \frac{\delta}{\alpha'}}\! \!\eta_r\d r\\
 &\leq 2K_2^2K_3\|\nabla f\|_\infty \frac{\delta}{\alpha'}\int_0^\delta\!(r+\sqrt{\veps' r})\d r+ \alpha' K_2\|\nabla f\|_{\infty}\theta(Y_{s_k}^{\veps',\alpha'})\int_0^\infty \eta_r\d r\\
 &\leq 2K_2^2K_3\|\nabla f\|_\infty \frac{\delta}{\alpha'}\Big(\frac{\delta^2}{2}+ \frac{2\sqrt{\veps' \delta^3}}{3}\Big)+\alpha' K_2\|\nabla f\|_\infty \theta(Y_{s_k}^{\veps',\alpha'})\int_0^\infty \!\eta_r\d r.
\end{aligned}
\end{equation}

Applying  condition \eqref{cond-1}   and It\^o's formula,
\begin{equation*}
\begin{split}
  \E\big[\theta(Y_t^{\veps',\alpha'})\big|\F_s\big]
  &=
  \theta(Y_s^{\veps',\alpha'})+\int_s^t\!\E \big[\frac1{\alpha'}Q(X_r^{ \veps',\alpha'})\theta(Y_r^{\veps',\alpha'}) \big|\F_s\big]\d r\\
  &\leq \theta(Y_s^{\veps',\alpha'})+\int_s^t\!\frac1{\alpha'} \big(-c_2\E\big[\theta(Y_r^{\veps',\alpha'})\big|\F_s\big] +c_3\big) \d r,\quad   t\geq s.
\end{split}
\end{equation*}
Then Gronwall's inequality yields
\begin{equation}\label{ab-9}
\E\big[\theta(Y_t^{\veps',\alpha'})\big|\F_s\big]
\leq \theta(Y_s^{\veps',\alpha'}) \e^{-c_2(t-s)/\alpha'}+\frac{c_3}{c_2},\quad  \,t\geq s.
\end{equation}

Inserting \eqref{ab-9}, \eqref{ab-8} into \eqref{ab-6.5}, we obtain that
\begin{equation*}
\begin{split}
\Upsilon_2&\leq   K_2^2K_3\|\nabla f\|_\infty t\frac{\delta}{\alpha'}\Big(\frac{\delta}{2} \!+\! \frac {2\sqrt{\veps'\delta}}{3}\Big)\!+\!\alpha'K_2\| \nabla f\|_\infty \! \!\int_0^\infty \!\!\eta_r\d r\!\sum_{k=0}^{N_t} \E\big[\theta(Y_{s_k}^{\veps', \alpha'})\big| \F_{s}\big]\\
&\leq K_2^2K_3\|\nabla f\|_\infty t\frac{\delta}{\alpha'}\Big(\frac{\delta}{2}\!+\! \frac {2\sqrt{\veps'\delta}}{3}\Big)\!\\
&\quad +\!\alpha'K_2\| \nabla f\|_\infty \!\int_0^\infty \!\!\eta_r\d r\Big(\theta(Y_s^{\veps',\alpha'}) \frac{1\!+\!c_2\delta}{c_2\delta}\!+\! \frac{c_3 (t\!+\!\delta)}{c_2\delta} \Big).
\end{split}
\end{equation*}
Taking $\delta=(\alpha')^{\frac 34}$, we arrive at
\begin{equation}\label{ab-10}
\begin{split}
  \Upsilon_2&\leq K_2^2K_3\|\nabla f\|_\infty t(\alpha')^{\frac 34}\Big(\frac{(\alpha')^{\frac 38}}{2}+\frac{2\sqrt{\veps'}}{3}\Big)\\
  &\quad + (\alpha')^{\frac 14}K_2\|\nabla f\|_\infty \int_0^\infty \!\eta_r\d r\!\Big(\theta(Y_s^{\veps',\alpha'})(1\!+\! c_2(\alpha')^{\frac 34})\!+\!\frac{c_3}{c_2}\big(t +(\alpha')^{\frac34}\big)\Big),
\end{split}
\end{equation}
which yields that $\lim_{\veps',\alpha'\to 0}\Upsilon_2=0$. Consequently, \eqref{ab-5} holds, and further $(\wt{X}_t)$ is a solution to ODE \eqref{a-4}.

If ODE \eqref{a-4} admits a unique solution, then for any subsequence of $(X_t^{\veps',\alpha'})_{t\in [0,T]}$ as $\veps',\alpha'\to0$, the tightness of $\{\mathcal{L}_{X^\ea};\veps,\alpha>0\}$ proved above tells us that there is further a subsequence of $(X_t^{\veps',\alpha'})_{t\in [0,T]}$,   which converges weakly to the unique solution $(\bar{X}_t)_{t\in [0,T]}$. Hence, the arbitrariness of the subsequence $(X_t^{\veps',\alpha'})_{t\in [0,T]}$ means that the whole sequence $(X_t^\ea)_{t\in [0,T]}$ converges weakly to $(\bar{X}_t)_{t\in [0,T]}$ as $\veps,\alpha\to 0$. The proof of Theorem \ref{thm-1.5} is complete.\fin

\end{document}